\documentclass[12pt]{amsart}
\usepackage{amscd,amssymb,amsmath}
\usepackage[all]{xy}
\usepackage[raiselinks=false,colorlinks=true,citecolor=blue,urlcolor=blue,
linkcolor=blue,bookmarksopen=true,pdftex]{hyperref}

\usepackage[mathscr]{eucal}
\newtheorem{theorem}{Theorem}[section]
\newtheorem{proposition}[theorem]{Proposition}
\newtheorem{lemma}[theorem]{Lemma}
\newtheorem{remark}[theorem]{Remark}
\newtheorem{definition}[theorem]{Definition}

\newcommand{\Mod}[1]{\ (\mathrm{mod}\ #1)}

\numberwithin{equation}{section}
\begin{document}

\title{K-theory of real Grassmann manifolds}
\author[S. Podder]{Sudeep Podder}
\address{Department of Mathematics, Indian Institute of Technology Madras, Chennai 600036, India}
\email{sudeep@smail.iitm.ac.in}
\author[P. Sankaran]{Parameswaran Sankaran}
\address{Chennai Mathematical Institute, SIPCOT IT Park, Siruseri, Kelambakkam, 603103, India}
\email{sankaran@cmi.ac.in}

\date{\today}

\newcommand{\af}{\alpha}
\newcommand{\et}{\eta}
\newcommand{\ga}{\gamma}
\newcommand{\ta}{\tau}
\newcommand{\ph}{\varphi}
\newcommand{\bt}{\beta}
\newcommand{\lb}{\lambda}
\newcommand{\wh}{\widehat}
\newcommand{\sg}{\sigma}
\newcommand{\om}{\omega}
\newcommand{\cH}{\mathcal H}
\newcommand{\cF}{\mathcal F}
\newcommand{\N}{\mathcal N}
\newcommand{\R}{\mathcal R}
\newcommand{\Ga}{\Gamma}
\newcommand{\cc}{\mathcal C}
\newcommand{\rank}{\mathrm{rank}}
\newcommand{\bea} {\begin{eqnarray*}}
\newcommand{\beq} {\begin{equation}}
\newcommand{\bey} {\begin{eqnarray}}
\newcommand{\eea} {\end{eqnarray*}}
\newcommand{\eeq} {\end{equation}}
\newcommand{\eey} {\end{eqnarray}}
\newcommand{\ovl}{\overline}
\newcommand{\vv}{\vspace{4mm}}
\newcommand{\lra}{\longrightarrow}
\newcommand{\SO}{\mathrm{SO}}
\renewcommand{\O}{\mathrm{O}}
\newcommand{\Spin}{\mathrm{Spin}}
\newcommand{\tor}{\mathrm{Tor}}

\begin{abstract}
Let $G_{n,k}$ denote the real Grassmann manifold of $k$-dimensional vector subspaces of $\mathbb R^n$.  
We compute the complex $K$-ring of $G_{n,k}$, up to a small indeterminacy,  
for all values of $n,k$ where $2\le k\le n-2$.  When $n\equiv 0\Mod {4}, k\equiv 1\Mod{2}$, we use the Hodgkin spectral sequence to determine the $K$-ring completely.   
\end{abstract}
\maketitle

\section{Introduction}
Let $G_{n,k}$ denote the real Grassmann manifold consisting of all $k$-dimensional vector subspaces in the real vector space $\mathbb R^n$.   We put the standard inner product on $\mathbb R^n$.  We have the identification of 
$G_{n,k}$ with the homogeneous space $\SO(n)/S(\mathrm{O}(k)\times \mathrm{O}(n-k)) $ where 
$\O(k)\times \O(n-k)$ is the subgroup of the orthogonal group $\O(n)$ that stabilises the subspace $\mathbb R^k$ 
spanned by the first $k$ standard basis vectors, and $S(\O(k)\times \O(n-k))=\SO(n)\cap (\O(k)\times \O(n-k))$. 
In this note our aim is to compute the complex $K$-ring of $G_{n,k}$.

Recall that the oriented 
Grassmann manifold $\widetilde{G}_{n,k}\cong \SO(n)/(\SO(k)\times \SO(n-k))$ is the double cover of $G_{n,k}$ and is simply-connected, except in the case of $\widetilde {G}_{2,1}\cong \mathbb S^1$.  The description of the $K$-ring of $\widetilde{G}_{n,k}$ goes back 
to work of Atiyah and Hirzebruch \cite{ah} when $n$ is odd or $k$ is even.   
 Note that in these cases, the subgroup 
$\SO(k)\times \SO(n-k)$ is connected and has rank equal to that of the whole group $\SO(n)$.  When $n$ is even and $k$ odd the $K$-ring was computed by Sankaran and Zvengrowski \cite{sz1}. 

The fact that $S(\O(k)\times \O(n-k))$ is not connected makes the determination of the ring $K(G_{n,k})$ difficult 
and, to the best of our knowledge, has not been carried out for $2\le k\le n-2$. 
Note that since $G_{n,k}\cong G_{n,n-k}$, it suffices to consider the case when $k\le n/2$. 
When $k=1$, $G_{n,1}$ is the same as the real projective space $\mathbb RP^{n-1},$ whose 
$K$-ring had been determined by Adams \cite{adams}. 

Our aim is to express $K^*(G_{n,k})=K^0(G_{n,k})\oplus K^1(G_{n,k})$ in terms of generators and relations.  
However, we have thus far only met with partial success.  We obtain complete results only under the assumption that 
$n\equiv 0\Mod{4}$ and 
$k$ odd.  In the remaining cases, our description is complete up to a small indeterminacy.   See Theorem \ref{main2} below 
and Proposition \ref{almostK}.    

We now state the two main results of this paper.  The proofs will be given in 
\S4 and \S5.

\begin{theorem}\label{main}
Let $n=2m, k=2s+1, n-k=2t+1$ and suppose that $m=s+t+1$ is even.  
Let $S$ denote the polynomial algebra $\mathbb Z[\lambda_1,\ldots,\lambda_s;\mu_1,\ldots,\mu_t,\theta] $ 
in $s+t+1$ variables.  Then
\[K^0(G_{n,k})=S/\mathcal I=\mathbb Z[\lambda_1,\ldots,\lambda_s;\mu_1,\ldots,\mu_t,\theta]/\mathcal I\]
where the ideal $\mathcal I$ is generated by the following elements:\\
(i)  $\theta^2-1,2^{m-1}(\theta-1),$\\
(ii)  $\sum_{0\le p\le j}\lambda_p\mu_{j-p}-{n\choose j}\theta^j, 1\le j\le m-1$, where $\lambda_{k-p}=\lambda_p,\mu_{n-k-q}=\mu_q$.\\
The $K^0(G_{n,k})$-module $K^1(G_{n,k})$ is the ideal generated by $\theta+1$ in the ring 
$S/\widetilde {\mathcal I}$, where $\widetilde{\mathcal  I}$ is generated by elements listed in (ii) above together with $\theta^2-1$. 
\end{theorem}

The element $[\theta]$ in the above theorem corresponds to the complexification of the Hopf line bundle $\xi=\xi_{n,k}$ over $G_{n,k},$ which is associated to double cover $\tilde G_{n,k}\to G_{n,k}$.    Note that since $\theta^2-1\in\widetilde{\mathcal I}$ we have $(\theta-1)\cdot y=0$ for all $y\in K^1(G_{n,k})$. It follows that the $S/\widetilde{\mathcal I}$-module $K^1(G_{n,k})$ is indeed a module over $S/\mathcal I=K^0(G_{n,k})$-module. 

Let $\gamma_{n,k}$ be the canonical (real) $k$-plane bundle over $G_{n,k}$.   Denote by $\mathcal K_{n,k}$ the 
$\lambda$-subring of $K(G_{n,k})$ generated by the class $[\gamma_{n,k}\otimes \mathbb C]$.  
An algebraic description $\mathcal K_{n,k}$ will be given in \S5.  

\begin{theorem}\label{main2}  Let $2\le k\le n/2$.  With the above notation, the inclusion $\mathcal K_{n,k}\hookrightarrow K(G_{n,k})$ has finite cokernel.
\end{theorem}

The main tool needed in the proof of Theorem \ref{main}  is the Hodgkin spectral sequence.  This will be recalled in \S2.
We need to compute the complex representation ring $RH_{n,k}$ of a certain subgroup $H_{n,k}$ of the spin group $\Spin(n)$ 
and determine its structure as a module over $R\Spin(n)$.  The relevant 
subgroup $H_{n,k}$ is such that $G_{n,k}\cong \Spin(n)/H_{n,k}$.   This is carried out in \S4 when $n\equiv 0\Mod{4}$ and $k$ is odd.  This seems rather complicated for arbitrary values of $n,k$.   As an application we obtain bounds for the order of the element $[\xi\otimes\mathbb C]-1\in K(G_{n,k})$ 
for any $n,k,~2\le k\le n/2$.  

Our proof of Theorem \ref{main2} uses standard arguments involving the Chern character.  

The Hodgkin spectral sequence had been used to determine the $K$-theory of 
many compact homogeneous manifolds.   Hodgkin \cite[\S12]{hodgkin} applied it 
to determine the $K$-ring of most of the compact simple Lie groups which are not necessarily simply connected.   Roux \cite{roux} used it to compute 
the $K$-ring of real Stiefel manifolds, independently of Gitler and Lam \cite{gitler-lam}, who had determined the same using a different approach.   
Antoniano, et al. \cite{aguz} and Barufatti and Hacon \cite{barufatti-hacon} used
the Hodgkin spectral sequence for computing the $K$-ring of real projective Stiefel manifolds, and  
Minami \cite{minami} 
for simply connected compact symmetric spaces.
See also \cite{sz1}, \cite{sz2}.

\noindent
{\bf Acknowledgments:} The authors thank the referee for his/her comments. 
The research of the first author was supported by University Grant Commission, India.
The research of the second author was partially supported by the Infosys Foundation.

\section{The Hodgkin spectral sequence}
We briefly recall the Hodgkin spectral sequence here.  Let $H$
 be a proper
closed subgroup of a compact Lie group $G$.   We denote the complex representation ring of 
$G$ by $RG$.   
Let $\rho:RG\to RH$ denote the restriction homomorphism and regard $RH$ as an $RG$-module via $\rho$.  
Hodgkin \cite{hodgkin} established the existence of a spectral sequence,  
whose $E_2$-diagram is given by $\tor^*_{RG}(RH,\mathbb Z)$, 
which converges to $K^*(G/H)$ when $\pi_1(G)$ is torsion-free.  Here $\tor^p_A(B,M)$ denotes $\tor_{-p}^A(B,M)$. In particular,  $\tor^*_{RG}(RH,\mathbb Z)$ is graded by non-positive integers.  
{\it We define the degree of an element $x\in \tor_A^p(B,M)$ to be $p$.}

 When the rings $RG, RH,$ and $\mathbb Z$ are given the trivial $\mathbb Z_2$ grading, we obtain a 
$\mathbb Z_2$-grading on $E_2^{p,q}$, where $E_2^{p,q}=\tor_{RG}^p(RH,\mathbb Z)$ if $q$ is even and 
is zero if $q$ is odd. In particular, $0=E_2^{p,q}=E_\infty^{p,q}$ if $q$ is odd.   The differential $d_r: E_r^{p,q}\to E_r^{p+r, q-r+1}$ vanishes when $r$ is even.  

Using the multiplication in $RH$, one obtains a $\mathbb Z_2$-graded ring structure 
on $\tor_{RG}^*(RH,\mathbb Z)$.  The differential in the spectral sequence is an anti-derivation, leading to a $\mathbb Z$-graded ring structure on $E_\infty^{*}$ which is compatible with the $\mathbb Z_2$-graded ring $K^*(G/H)$.

 If $\tor^*_{RG}(RH,\mathbb Z)$ is generated by elements of degree at least $-2$, then the spectral sequence collapses at the $E_2$-stage and we have $\tor^*_{RG}(RH,\mathbb Z)\cong K^*(G/H)$.  See \cite{roux}.

Pittie \cite{pittie} has shown that $RH$ is stably free over $RG$ if $H$ is  connected, $\pi_1(G)$ is torsion-free, and the 
rank of $H$ equals the rank of $G$, i.e., if $H$ has a maximal torus $T\subset H$ which is maximal in $G$.  Moreover, if 
$|W(G,T)|/|W(H,T)|>1+\dim T$, then $RH$ is a free $RG$-module.  (Here $W(G,T)$ denotes the Weyl group of $G$ with respect to $T$.)  Consequently the Hodgkin spectral sequence collapses and we have 
$K(G/H)=\tor^0_{RG}(RH;\mathbb Z)=RH\otimes_{RG}\mathbb Z$.  In case $G$ is prime to 
the exceptional Lie groups of type $E_6, E_7, E_8,$ this was proved by Atiyah and Hirzebruch \cite{ah}, who  conjectured 
its validity for any $G$ with $\pi_1(G)$ torsion-free.   

\subsection{Change of rings spectral sequence}\label{crss}
Suppose that $G$ is simply connected so that $RG$ is a polynomial ring $\mathbb Z[x_1,\ldots, x_m]$.  
When $RH$ is not a free $RG$-module (via the restriction homomorphism),  but is free over a subring 
$\Lambda =\mathbb Z[x_1,\ldots, x_r]$, then it is possible to use the change of rings spectral sequence 
due to Cartan and Eilenberg \cite{cartan-eilenberg} to compute $\tor^*_{RG}(RH,\mathbb Z)$.   
See \cite{roux}, \cite[\S6]{aguz} and also \cite[\S6]{barufatti-hacon} for a more 
detailed discussion on the use of the change of rings spectral sequence in the computation of $K(G/H)$.
We now recall the Cartan-Eilenberg change of rings theorem.

Let $K$ be any ring.  
 A $K$-algebra $\Lambda$ together with a $K$-homomorphism $\varepsilon:\Lambda\to K$ is called a {\it supplemented $K$-algebra} with augmentation $\varepsilon$.  Let $(\Lambda,\varepsilon), (\Gamma,\eta)$ be supplemented $K$-algebras, and let 
$\varphi:\Lambda\to \Gamma$ be a $K$-algebra homomorphism such that $\eta\circ \varphi=\varepsilon$.  
Denote $\ker(\varepsilon)$ by $I(\Lambda)$.  
 A $K$-algebra homomorphism $\varphi:\Lambda
\to \Gamma$ is {\it normal} if the left ideal, denoted $\Gamma\cdot I(\Lambda)$, of 
$\Gamma$ generated by $\varphi(I(\Lambda)),$ is also a 
right ideal of $\Gamma$ (always the case when $K$ is commutative).
Then $\Omega:=\Gamma/(\Gamma\cdot I(\Lambda))$ is a supplemented $K$-algebra.

We shall state the theorem in the special case of {\it commutative} augmented $K$-algebras. 
 So if $\Gamma, ~\Lambda$ are supplemented, any augmentation preserving $K$-homomorphism $\Gamma\to \Lambda$ is normal.  
  In our applications, $K=\mathbb Z$, $\Gamma=RG$, $~\Lambda$ will be a subring of $\Gamma$, and $A=RH$, where the $\Gamma$-module structure is given 
via the restriction homomorphism $\rho:RG\to RH$.   Also, the $\Omega$-module 
$C$ in the statement of the theorem below will be $\mathbb Z$ (via the augmentation).

\begin{theorem} \label{changeofrings} {\em (\cite[Theorem 6.1, Chapter XVI]{cartan-eilenberg})}  We keep the above notations. Suppose that $K$ is commutative. 
Suppose that $\varphi: \Lambda\to \Gamma$ is normal and that $\Gamma$ is projective as a $\Lambda$-module (via $\varphi)$.  Then, for any $\Gamma$-module $A$ and $\Omega$-module $C$, there exists a spectral sequence 
$\tor^\Omega_*(\mathrm{Tor}^\Lambda_*(A,K),C)$ that converges to $\mathrm{Tor}_*^\Gamma(A,C)$.  
\end{theorem}

The $\Omega$-module structure on $\textrm{Tor}^\Lambda_q(A,K)$ arises from the functorial isomorphism 
$\textrm{Tor}^\Gamma_q(A,\Omega)=\textrm{Tor}^\Gamma_q(A,\Gamma\otimes_\Lambda K)\cong 
\textrm{Tor}^\Lambda_q(A,K)$. (See \cite{cartan-eilenberg} for details.)


\section{The representation ring of $H_{n,k}$}
We follow the notations of Husemoller's book \cite{husemoller} closely in our 
description of the representation rings of the groups $\SO(n)$ and $\Spin(n)$.  

Let $2\le k\le \lfloor n/2\rfloor$.
Recall that $H_{n,k}$ is the inverse image of $S(\O(k)\times \O(n-k))$ under the double cover $\pi: \Spin(n)\to \SO(n)$.
The identity component of $H_{n,k}$ is the group $H^0_{n,k}:=\Spin(k)\cdot \Spin(n-k)\subset \Spin(n)$ with quotient 
$H_{n,k}/H_{n,k}^0\cong \mathbb Z_2$.   Although the representation ring of $H^0_{n,k}$ has been worked out in 
\cite{sz1}, we shall give most of the details here in order to make the exposition self-contained. 
Note that $H^0_{n,k}$ is the quotient of $\Spin(k)\times \Spin(n-k) $ 
by the cyclic subgroup of order $2$ generated by $(-1,-1)$.    The canonical surjection $\Spin(k)\times \Spin(n-k)\to H^0_{n,k}$ induces a ring 
monomorphism $RH^0_{n,k}\to R(\Spin(k)\times \Spin(n-k))$ which we regard as an inclusion.  
The image is generated as an abelian group by representations of $\Spin(k)\times \Spin(n-k)$ on which $(-1,-1)$ acts as identity.  Likewise, the projection 
$H^0_{n,k}\to \SO(k)\times \SO(n-k)$ induces a monomorphism $R(\SO(k)\times \SO(n-k))\to RH^0_{n,k}$, which we regard as an inclusion, whose 
image is generated by representations of $H^0_{n,k}$ on which the kernel of the projection acts as the identity.  This allows us to describe $RH^0_{n,k}$ in a straightforward manner.   The ring $R(\SO(k)\times \SO(n-k))$ is a polynomial ring when $n$ is even and $k$ is odd.  The ring homomorphism 
$R\SO(2r+1)\to R\SO(2r)$ induced by the inclusion $\SO(2r)\hookrightarrow \SO(2r+1)$ is a monomorphism.  Moreover, $R\SO(2r+1)$ is a polynomial 
ring in $r$ indeterminates.  The ring $R\SO(2r)$ is not isomorphic to a polynomial algebra; it is known that 
$R\SO(2r)$ is generated over $R\SO(2r+1)$ by an element $\lambda^+_r$ which satisfies a monic quadratic 
equation.  As such $R\SO(2r)$ is a free 
$R\SO(2r+1)$-module  of rank $2$.    So,  for all parities of $k,n$, $R(\SO(k)\times \SO(n-k))$ is a free module of finite rank over a polynomial 
ring generated by $\lfloor k/2\rfloor+\lfloor (n-k)/2\rfloor$ indeterminates.  
We will show in this section that the same statement holds for $RH_{n,k}$ as well.

Before proceeding further in describing $RH^0_{n,k}, RH_{n,k},$ we need to introduce notations for 
certain natural representations of the spin and special orthogonal groups.

Set $k=2s+\varepsilon, n-k=2t+\eta, \varepsilon, \eta\in \{0,1\}$ where $s,t$ are integers.  Now $n= 2s+2t+1$ if $n$ is odd. When $n$ is even, both $k$ and $ n-k$ are of same parity and 
$n=2s+2t$ or $n=2s+2t+2$ according as $k$ is even or odd. 
Let $\lambda_1$ denote the standard $k$-dimensional complex representation of $\SO(k)$.  We denote by $\lambda_j\in R\SO(k)$ the $j$th exterior 
power $\Lambda^j_\mathbb C(\lambda_1),  j\le k$.  (It is understood that $\lambda_0=1$, the trivial representation).  \footnote{We shall often use the same notation for a representation and its class in the representation ring.}
We have the equality
\[\lambda_j=\lambda_{k-j}~\mathrm{~in~} R\SO(k).\]  

When $k$ is even, the Hodge star operator $*$ yields a splitting $\lambda_s=\lambda^+_s+ \lambda^-_s,$ 
where $\lambda^+_s,\lambda^-_s\in R\SO(2s)$ are the classes of $+1, -1$-eigenspaces when $k\equiv 0\pmod{4}$ and are the $i, -i$-eigenspaces when $k\equiv 2 \Mod{4}$ respectively. 
In the case of $\Spin(k)$ we have the spin representation $\Delta_s$.  When $k$ is even, it splits as a sum of two half-spin representations $\Delta^+_s,\Delta_s^-$; they are distinguished by the way an element $z_0$ in the centre of $\Spin(k)$ acts.  (This will be made precise later.)  We have the following 
theorem proved in \cite[\S10, Chapter 13]{husemoller}.  In the case of 
$R\SO(2s)$, our description is slightly different from the one given in Husemoller's book {\it op. cit.}, but it is readily seen that the two descriptions are equivalent. 

\begin{theorem} \label{gen-RG}
With the above notations,  we have \\
(i)  $R\Spin(2s)=\mathbb Z[\lambda_1,\cdots, \lambda_{s-2}, \Delta^+_s,\Delta^-_s]$, \\
(ii) $R\Spin(2s+1)=\mathbb Z[\lambda_1,\ldots, \lambda_{s-1}, \Delta_s],$ \\
(iii) $ R\SO(2s+1)=\mathbb Z[\lambda_1,\ldots,\lambda_s],$ and,  \\
(iv) $R\SO(2s)=\mathbb Z[\lambda_1,\lambda_2,\ldots, \lambda_{s}][\lambda_s^+]/\!\!\sim$ where the ideal of relations is generated by 
$ (\lambda^{+}_s)^2-a\lambda^+_s-b$ for 
suitable polynomials $a,b$ in $\lambda_j, 1\le j\le s$ (with $\mathbb Z$-coefficients). 
\end{theorem}
As the notation suggests, the rings $R\Spin(2s), R\Spin(2s+1), R\SO(2s+1)$ are polynomial rings in the indicated variables.  Also, the elements $\lambda_j, 1\le j\le s,$ in $R\SO(2s)$ are algebraically independent.   

\begin{remark}
The quadratic relation 
that 
$\lambda_s^+$ 
satisfies over $\mathbb Z[\lambda_1,\ldots,\lambda_s]$ can be explicitly 
written down as follows:  Set $\lambda^-_s:=\lambda_s-\lambda^+_s$.  
From \cite[Theorem 10.3, Chapter 13]{husemoller}, we have the relation 
$\lambda_s^+\cdot \lambda^-_s=(\lambda_{s-1}+\lambda_{s-3}+\cdots)^2 
-\lambda_s(\lambda_{s-2}+\lambda_{s-4}+\cdots)-(\lambda_{s-2}+\lambda_{s_4}+\cdots)^2\in \mathbb Z[\lambda_1,\ldots, \lambda_s]$.   Denoting the {\em  negative} of the 
right hand side of the last equality by $b$ and setting 
$a:=\lambda_s$ we have 
\[ (\lambda^+_s)^2=\lambda^+_s(\lambda_s-\lambda^-_s)=a\lambda^+_s+b.\]
\end{remark}

The inclusion $\Spin(2s)\hookrightarrow \Spin(2s+1)$ induces an injective ring homomorphism $\rho: R\Spin(2s+1)\to R\Spin(2s)$ where $\rho(\Delta_{s})=\Delta_s^++\Delta_s^-, \rho(\lambda_i)=\lambda_i+\lambda_{i-1}, 
1\le i\le s.$
The homomorphism $R\Spin(2s)\to R\Spin(2s-1)$ induced by the inclusion $\Spin(2s-1)\hookrightarrow \Spin(2s)$ 
is given by $\lambda_j\mapsto \lambda_j+\lambda_{j-1}, 1\le j<s, \Delta^\pm_s\mapsto 
\Delta_{s-1}$.   These restriction homomorphisms also yield the restrictions $R\SO(k)\to R\SO(k-1)$ for any parity of $k$.  

Recall that given any two compact Lie groups $H_1, H_2$, we have $R (H_1\times H_2)=RH_1\otimes RH_2$.  We have the  natural quotient homomorphisms $\pi_0: \Spin(k)\times\Spin(n-k)\to H^0_{n,k}$ and $\pi:H^0_{n,k}\to \SO(k)\times \SO(n-k)$ where $\ker(\pi_0)\cong \mathbb Z_2$ is generated by $(-1,-1)\in \Spin(k)\times \Spin(n-k)$ and 
$\ker \pi\cong \mathbb Z_2$, by $\pi_0(1,-1)\in H^0_{n,k}$.   
We shall regard the ring homomorphisms $\pi_0^*:RH^0_{n,k}\to R(\Spin(k)\times 
\Spin(n-k)),  ~\pi^*:R(\SO(k)\times \SO(n-k))\to RH^0_{n,k}$, which are injective, as inclusions.  
It is easy to see that $RH^0_{n,k}$ is generated as an 
$R(\SO(k)\times \SO(n-k))$-algebra by elements $xy\in R(\Spin(k)\times \Spin(n-k))$ where $x,y$ vary over the $R(\SO(k)\times \SO(n-k))$-algebra generators of $R(\Spin(k)\times \Spin(n-k))$.   The following description, in Proposition \ref{rh0nk}, of $RH^0_{n,k}$ is an immediate consequence of Theorem \ref{gen-RG}.

We shall  use the notation $\mu_j\in R\SO(n-k)$ for the element represented by the $j$th exterior power of the standard representation 
of $\SO(n-k)$.  Also $\Delta'_t, $ and $\Delta^{'\pm}_t$ will denote the spin and half-spin representations of $\Spin(n-k)$ respectively.  Thus $R(\SO(k)\times 
\SO(n-k))$ contains the polynomial subring $\mathbb Z[\lambda_1, \ldots, \lambda_s, \mu_1,\ldots, \mu_t].$

\begin{proposition} \label{rh0nk}  We keep the above notations.  
Let $R:=R(\SO(k)\times \SO(n-k))$.  Then
\[ RH^0_{n,k}=\begin{cases}
R[\Delta_s\Delta_t'] ,& \mbox{if~} k=2s+1, n-k=2t+1,\\
R[\Delta_s(\Delta'_t)^{\pm}], & \mbox{if~} k=2s+1, n-k=2t,\\
R[\Delta_s^\pm \Delta'_t], & \mbox{\if~} k=2s, n-k=2t+1,\\ 
R[\Delta_s^\pm (\Delta_t')^\pm, \Delta^\pm_s(\Delta'_t)^\mp], &\mbox{if~} k=2s,n-k=2t.\\
\end{cases} \]
Moreover, the squares of the indicated generators belong to $R$. 
\end{proposition}

\noindent
{\bf Notations:}  We shall denote by $\Delta_{s,t}$ the element $\Delta_s\Delta'_t\in R(\Spin(k)\times \Spin(n-k))$.  Also $\Delta_{s,t}^{\varepsilon,\eta}$ will denote 
$\Delta_s^\varepsilon\cdot(\Delta_t')^\eta, \varepsilon,\eta\in \{+,-\}.$  Also, we shall use upper case letters $\Lambda_j, 1\le j\le m$, etc.,  to denote the generators of $R\Spin(n)$ and  similarly $\lambda_1,\ldots, \lambda_s$ (resp. $\mu_1,\ldots, \mu_t$) 
 to denote generators of $R\Spin(k)$ (resp. $R\Spin(n-k)$) as in Theorem \ref{gen-RG}.

Next we turn our attention to  the representation ring of $H_{n,k}$.   Recall 
that $2\le k\le n/2$ and so $n\ge 4$.
First we analyse when the exact sequence 
\[ 1\to H^0_{n,k}\to H_{n,k}\to  Z\to 1 \eqno(1)\]
splits.    Evidently, the sequence splits if and only if there exists an element $z_0\in H_{n,k} \subset \Spin(n)$ of order $2$ such that $z_0\notin H^0_{n,k}$.  
Taking  $z_0:=e_1e_2e_3e_n\in C_n$, we see that $z_0^2=1$ and $z_0\in H_{n,k}\setminus H^0_{n,k}$, so the short exact sequence (1) splits.  
  Here $C_n$ denotes the Clifford algebra of the quadratic space 
$(\mathbb R^n, -||\cdot||^2)$ and $e_1,\ldots,e_n$ denote the standard basis vectors of $\mathbb R^n$.  So $H_{n,k}\cong H^0_{n,k}\rtimes \mathbb Z_2.$  

Suppose that 
$H_{n,k}=H^0_{n,k}\times Z$ and let $z_0$ be the generator of 
$Z\cong \mathbb Z_2$.   
Then $\pi(H_{n,k}^0)\times \pi(Z)=\pi(H_{n,k})=S(\textrm{O}(k)\times \textrm{O}(n-k))$ is isomorphic to the product 
$\SO(k)\times \SO(n-k)\times \{\pm I_n\} $.  In particular, $n$ is even and $k$ is odd and $z_0\in Z$ maps to $-I_n$.  So,  the order $2$ element $z_0$ is in the centre of $\Spin(n)$.   It follows that $n\equiv 0\Mod 4, k\equiv 1 \Mod 2$.

When $n\equiv 0\Mod{4},~k\equiv 1\Mod{2}$, we may take $z_0=e_1e_2\cdots e_n\in H_{n,k}$.  Then $z_0$ is in the centre of $H_{n,k}$ and  $z_0\notin H_{n,k}^0$ and so $H_{n,k}$ is the direct product $H_{n,k}^0\times Z$.   

Thus $H_{n,k}\cong H_{n,k}^0\times \mathbb Z_2$ if and only if $n\equiv 0\Mod4, k\equiv 1\Mod2.$
 
Using Proposition \ref{rh0nk}, we obtain the following.  

\begin{proposition} \label{rhnk}  We keep the above notations.   Let $k=2s+1,n-k=2t+1, s+t $ odd.  Let $f_{s,t}\in R:=R(\SO(k)\times \SO(n-k))$ be the element such that $\Delta_{s,t}^2=f_{s,t}$ and let $\theta$ be the class of the unique non-trivial one-dimensional representation of $H_{s,t}$.  Then 
\[RH_{n,k}=RH^0_{n,k}\otimes RZ=R[\Delta_{s,t},\theta]/\langle \theta^2-1,\Delta_{s,t}^2-f_{s,t}\rangle. \eqno(2)\]
 In particular, $RH_{n,k}$ is a free $R$-module 
with basis $\{1,\theta ,\Delta_{s,t},\theta\Delta_{s,t}\}.$
\end{proposition}

Writing $\lambda_0=1=\mu_0$,  $f_{s,t}\in R$ can be expressed as a polynomial in $\lambda_p, \mu_q\in R, 0\le p\le s, 0\le q\le t$ as follows (see \cite[Theorem 10.3, Chapter 14]{husemoller}.)
\[f_{s,t}=\Delta_{s,t}^2=\Delta_s^2\cdot(\Delta_t')^{2}=(\sum_{0\le p\le s}\lambda_p)(\sum_{0\le q\le t}\mu_q)=
\sum_{0\le r\le s+t}(\sum_{p+q=r} \lambda_p\mu_q). \eqno(3)\]

 
\section{The restriction homomorphism $R\Spin(n)\to RH_{n,k}$}
Throughout this section we assume that $k=2s+1, n-k=2t+1$ so that $n=2m,$ where $m:=s+t+1$.  Also we shall assume that $s+t$ is odd so that $n\equiv 0\Mod{4}$.   Hence $H_{n,k}=H_{n,k}^0\times Z$ where $Z\cong \mathbb Z_2$ is generated by $z_0=e_1\cdots e_n\in \Spin(n)$ .

The double covering $\phi: \Spin(n) \to \SO(n)$ is defined 
as $\phi(u)(x)=uxu^*,~x\in \mathbb R^n,$ where $*$ is (the restriction to $\Spin(n)$ of) the anti-involution of the Clifford algebra $C_n,$ uniquely defined by the requirement:  $v^*=v,~v\in \mathbb R^n$.  We refer 
the reader to \cite{husemoller} concerning the spin group and its representation ring.

\subsection*{Maximal tori}
Set $\omega(\theta_1,\ldots, \theta_m):= \prod_{1\le j\le m}(\cos2\pi  \theta_j+\sin 2\pi \theta_j . e_{2j-1}e_{2j}) \in \Spin(n)$ where $ \theta_j\in \mathbb R$.  
Then $\widetilde{T}:=\{\omega(\theta_1,\ldots, \theta_m)\in \Spin(n)\mid \theta_j\in \mathbb R, 1\le j\le m\}\cong (\mathbb S^1)^m$ is a maximal torus of $\Spin(n)$.  Its image in $\SO(n)$ is the standard maximal torus 
$T:=\SO(2)\times \cdots \times \SO(2)$ whose elements restrict to rotations on $\mathbb R e_{2j-1}+\mathbb R e_{2j}, 1\le j\le m.$ 
 In fact $\phi(\omega(\theta_1,\ldots, \theta_m))=
D(2\theta_1,\ldots, 2\theta_m)\in T$ where $D(t_1,\ldots, t_m)$ restricts to the positive rotation by angle $2\pi t_j$ on the  oriented vector subspace $\mathbb Re_{2j-1}+\mathbb R e_{2j}, 1\le j\le m$, the orientation being given by the 
ordering $e_{2j-1},e_{2j}$ of the basis elements.

Let $\mathbb T$ be 
the `standard torus' $(\mathbb S^1)^m =(\mathbb R/\mathbb Z) ^m$.  One has a homomorphism $\omega: \mathbb T\to \widetilde T$ defined by $(\theta_1,\ldots,\theta_m)\mapsto \omega(\theta_1,\ldots, \theta_m)$.    Note that $\omega(\theta_1+\varepsilon_1/2, \ldots, \theta_m+\varepsilon_m/2)=(-1)^\varepsilon \omega(\theta_1,\ldots,\theta_m)$
where $\varepsilon_j\in \{0,1\}$ for all $ j,$ and $\varepsilon=\sum_{1\le j\le m}\varepsilon_j$.  In particular $\ker(\omega)\cong 
(\mathbb Z_2)^{m-1}$.  The kernel of $\phi\circ \omega:\mathbb T\to T$ is readily seen to be $\mathbb Z_2^m\cong \{-1,1\}^m\subset \mathbb T$.

Since 
$n$ is even and $k$ is odd, the rank of $H^0_{n,k}$ equals $m-1=\textrm{rank}(\Spin(n))-1$.  
 In this case, $\widetilde T_0:=H^0_{n,k}\cap \widetilde T=\{\omega(\theta_1,\ldots, \theta_m)\in \widetilde T\mid \theta_{s+1}=0\}$ is a maximal torus of $H^0_{n,k}$.
Also,  we observe that the element $z_0=e_1\ldots e_n,$ the generator of $Z$, belongs to $\widetilde T $.    Let 
$T_0=\pi(\widetilde T_0)= T\cap (\SO(k)\times \SO(n-k))$ which is a maximal torus of $\SO(k)\times \SO(n-k)$.

The representation rings of $\widetilde T, \widetilde T_0, T, T_0 $ are viewed as subrings of $R\mathbb T$ as follows:   Let $u_j:\mathbb T\to \mathbb S^1$ be the $j$th projection, regarded as a character.  We also denote the corresponding $1$-dimensional representation of $\mathbb T$ by the same symbol $u_j$.  Then 
$R\mathbb T=\mathbb Z[u_1^{\pm 1},\ldots, u_m^{\pm 1}]$, $R\widetilde T=\mathbb Z[u_1^{\pm 2},\ldots, u_m^{\pm 2},u_1\cdots u_m ],$ and, $
 RT=\mathbb Z[u_1^{\pm 2},\ldots, u_m^{\pm 2}]$, both regarded as subrings of $R\mathbb T$.  
 Also $H_{n,k}\cap \widetilde T=\widetilde T_0\times Z$.  We have 
 \[RT_0=\mathbb Z[u_1^{\pm 2},\ldots, u_s^{\pm 2}, v_1^{\pm 2},\ldots, v_{t}^{\pm 2}]\subset R\mathbb T\]
where $v_j:=u_{s+j+1}, 1\le j\le t, $ and,  
\[R\widetilde T_0=\mathbb Z[u_1^{\pm 2},\ldots, u_s^{\pm 2}, v_1^{\pm 2},\ldots,v_t^{\pm 2},  u_1\cdots u_s v_1\cdots v_t]
\subset R\widetilde T. \]

In order to determine the restriction homomorphism $\rho: R\Spin(n)\to RH_{n,k}$, we first consider the 
homomorphism $R\Spin(n)\to R\Spin (n)\otimes RZ$ induced by the homomorphism 
$\mu:  \Spin(n)\times Z\to \Spin(n)$ defined by multiplication: $(g,z)\mapsto gz$.     
Note that the restriction 
of $\mu$ to $H^0_{n,k}\times Z$ is an isomorphism $H^0_{n,k}\times 
Z\to H_{n,k}$.
The homomorphisms $H^0_{n,k}\times Z\to H_{n,k}$, 
$\widetilde T\times Z\to \widetilde T, ~\widetilde{T}_0\times Z\to \widetilde T$, and $\widetilde T_0\times Z\to H_{n,k}$, each of which is obtained from $\mu$ by appropriately restricting its domain and co-domain, will all be denoted 
by the same symbol $\mu$ by an abuse of notation. 
These group homomorphisms induce homomorphisms of rings 
 $\mu^*: R\widetilde T 
\to R\widetilde T \otimes  RZ,~\mu^*:R\widetilde T\to R\widetilde T_0\otimes RZ$, $\mu^*:RH_{n,k}\to R\widetilde T_0\otimes RZ$, 
$ \mu^*: RH_{n,k}\stackrel{\cong}{\longrightarrow} RH^0_{n,k}\otimes RZ,$ and 
$\mu^*:R\Spin(n)\to R\Spin(n)\otimes RZ$.

Let $\sigma:\widetilde T\hookrightarrow \Spin(n)$ be the inclusion. 
We have the following commutative diagram where the homomorphisms in the first row are induced by respective inclusions of groups.  

\[\begin{array}{rcccccl}
  R\Spin(n) \otimes RZ &\hookrightarrow & R\widetilde T\otimes RZ &\rightarrow &R\widetilde T_0\otimes RZ\hookleftarrow&RH^0_{n,k}\otimes RZ&\\
\uparrow \mu^*~~~~~& & \mu^*\uparrow &&\uparrow  id&\uparrow \mu^*&\\
R\Spin(n) &\stackrel{\sigma^*}{\hookrightarrow }&R\widetilde T&\stackrel{\mu^*}{\longrightarrow} &R\widetilde T_0\otimes R Z \stackrel{~~\mu^*}{~~\hookleftarrow} &RH_{n,k}
\\
\end{array}\eqno(4)
\]

The inclusion $\sigma^*:R\Spin(n)\hookrightarrow R\widetilde T$ is via the identification of $R\Spin(n)$ with the invariant 
subgroup of $R\widetilde T$ under the action of the Weyl group $W(\Spin(n),\widetilde T)$.  
Similarly we have the inclusion $RH^0_{n,k}\hookrightarrow R\widetilde T_0$ which in turn induces $RH_{n,k}\hookrightarrow R\widetilde T_0\otimes RZ$. 
Moreover, 
$\mu^*(R\Spin(n))$ is contained in $RH_{n,k}\subset R\widetilde T_0\otimes RZ$ since $H_{n,k}\subset \Spin(n)$.  
This allows one to describe the restriction homomorphism 
$\rho: R\Spin(n) \to RH_{n,k}$ easily, once $\mu^*:R\widetilde T 
\to R\widetilde T_0\otimes RZ$ is determined.  This we shall carry out below, 
with $\theta$ as in Proposition \ref{rhnk}.

Routine computation, using $n=2m, m$ even, yields that 
\[u_1\cdots u_m(z_0)=\begin{cases} 1& \mbox{~if~} n\equiv 0\Mod{8}\\ \theta(z_0)&\mbox{~if~} n\equiv 4\Mod{8}.\end{cases}\eqno(5)\] 

When $t\in \widetilde T_0$, we have $u_{s+1}^2(t)=1$  
and so $u_1.\ldots.u_m$ restrict to $u_1\cdots u_s\cdot v_{1}\cdots v_t$ on $\widetilde T_0$.  Therefore 
\[\mu^*(u^{\pm 2}_j)=\begin{cases}
\theta u_j^{\pm 2}, & 1\le j\le s,\\
\theta, & j=s+1,\\
\theta v_{j-s-1}^{\pm 2}, &s+1<j\le m,\\
\end{cases}
\eqno (6)\]  
and,
\[
 \mu^* (u_1\cdots u_m)=\begin{cases} \prod _{1\le j\le s} u_j\cdot \prod_{1\le j\le t} v_j, & n\equiv 0\Mod{8},\\
\theta \prod_{1\le j\le s} u_j\cdot\prod _{1\le j\le t}v_j , &n\equiv 4\Mod{8}.
\end{cases}\eqno(7)\]

Let $e_j(x_1,\ldots, x_r)$ denote the $j$th elementary symmetric polynomial in $x_1,\ldots, x_r$.  Recall that 
$\sigma^*(\Lambda_j)=e_j(u_1^2,u_1^{-2}, \ldots, u_m^2, u_m^{-2})$.  So, for $1\le j\le m$, we have 
\[ \begin{array}{rcl}\rho(\Lambda_j)&=&\mu^*(e_j(u_1^2,
u_1^{-2},\ldots, u_m^2,u_m^{-2}))\\&=&\theta^je_j(u_1^2,u_1^{-2},\ldots, u_s^2, u_s^{-2}, 1,1,v_1^2,v_1^{-2},\ldots, v_t^2,v_t^{-2})\\
&=& \theta^j\sum_{p+q=j}e_p(u_1^2,u_1^{-2}, \ldots, u_s^2,u_s^{-2},1)\cdot e_q(v_1^2,v_1^{-2},\ldots, v_t^2,v_t^{-2}, 1)\\
&=&\theta^j \cdot \sum_{p+q=j;0\le p\le k,0\le q\le n-k}\lambda_p\mu_q,\\
&=&\theta^jf_j
\end{array}\eqno(8)\]  
for a suitable element $f_j=f_j(\lambda_1,\ldots, \lambda_s,\mu_1,\ldots,\mu_t)\in R$ 
since $\lambda_p=\lambda_{k-p}, \mu_q=\mu_{n-k-q}.$

Using Equations (6) and (7) we obtain that if $\varepsilon_j\in \{1,-1\}$, then 
\[\mu^*(u_1^{\varepsilon_1}\cdots u_m^{\varepsilon_m})
=\theta^\varepsilon  u_1^{\varepsilon_1}\cdots u_s^{\varepsilon_s}\cdot v_1^{\eta_1}\cdots v_t^{\eta_t}\eqno(9)\] where 
$ \eta_j=\varepsilon_{s+1+j}, $ and the value of $\varepsilon \in \{0,1\}$ is obtained as follows: $\varepsilon\equiv \sum_{1\le j\le m}\varepsilon_j\Mod{2}$ if $n\equiv 0\Mod{8}$ and $\varepsilon\equiv 1+ \sum_{1\le j\le m} \varepsilon_j\Mod{2}$ if 
$n\equiv 4\Mod{8}$.
  
The following proposition now follows immediately from equations (8), (9), and the definitions of $\Delta^\pm_m, ~\Delta_{s,t}$.  

\begin{proposition}  \label{restriction}
Let $n=2m\equiv 0\Mod{4}, k=2s+1, n-k=2t+1$.  
With the above notations, the restriction homomorphism 
$\rho:R\Spin(n)\to RH_{n,k}$ is defined by $\rho(\Lambda_j)=\Lambda'_j=\theta ^j\sum_{p+q=j} \lambda_p\mu_q=\theta^jf_j, 1\le j\le m-1$, 
$\rho(\Delta^+  _{m})
=\theta ^\varepsilon \Delta_{s,t}, \rho (\Delta^-_m)=\theta^{1+\varepsilon}\Delta_{s,t}$ where 
$\varepsilon =0, 1$ according as $n\equiv 0\Mod{8}$ or $n\equiv 4\Mod{ 8}$ respectively.
\end{proposition}

The ring 
$R':=\mathbb Z[\theta^p\lambda_p, \theta^q\mu_q; 1\le p\le s, 1\le q\le t] \subset RH_{n,k}$ is mapped to the polynomial 
ring $\mathbb Z[\lambda_p,\mu_q;1\le p\le s, 1\le q\le t]=R=R(\SO(k)\times \SO(n-k))$ by an automorphism of the ring $R[\theta]$ 
since $\theta$ is invertible.  It follows that $R'$ is a polynomial ring in $s+t=m-1$ indeterminates.   Evidently, $R'[\theta]=R[\theta]$.  

\begin{lemma} \label{freenessoverLambda}  Let $n=2m\equiv 0\Mod{4}, k=2s+1, n-k=2t+1$.  
Let  $R'[\theta]=R[\theta]\subset RH_{n,k}.$   Then $R'[\theta]$ 
is a free $\Lambda'$-module of rank $2{m-1\choose s}$ where
$\Lambda':=\mathbb Z[\Lambda_1',\ldots, \Lambda_{m-1}']\subset RH_{n,k}$.  
 In particular, 
$\Lambda'_1,\ldots, \Lambda_{m-1}'$ are algebraically independent. 
Also 
$RH_{n,k}=R[\theta,\Delta_{s,t}]$ is a free module of rank $4{m-1\choose s}$ over $\mathbb Z[\Lambda_1,\ldots, \Lambda_{m-1}]$ via $\rho$.   
\end{lemma}
\begin{proof}
Since $R=R(\SO(k)\times \SO(n-k))$ is a polynomial algebra in $s+t=m-1$ indeterminates, the algebraic independence of $\Lambda_1',\ldots, \Lambda_{m-1}'$ would follow once we show that $R[\theta]\cong R\oplus R$ 
is a finitely generated free $\Lambda'$-module. 

First note that $\Lambda'[\theta]$ is free as a $\Lambda'$-module with basis $\{1,\theta\}$.  

Next we will show that $R[\theta]\subset RH_{n,k}$ is free as a $\Lambda'[\theta]$-module of rank ${m-1\choose s}$.  Let $\Lambda_0=\mathbb Z[f_1,\ldots, f_{m-1}]$. 
Then $\Lambda'[\theta]=\Lambda_0[\theta]=\Lambda_0\otimes_\mathbb Z\mathbb Z[\theta]$. 
Since $R[\theta]=R\otimes_\mathbb Z
\mathbb Z[\theta],$ it suffices to show that $R$ is free as a module over $\Lambda_0\subset R$, of 
rank ${m-1\choose s}$.

Denote by $\rho_0: R\Spin(n)\to RH_{n,k}\to RH^0_{n,k}$ the restriction homomorphism induced by 
the inclusion $H_{n,k}^0\hookrightarrow H_{n,k}\hookrightarrow \Spin(n)$.    Then $\Lambda_0=\rho_0(\Lambda)$ and $ \rho_0(\Lambda)\subset R\subset R[\Delta_{s,t}]$.  Then $R$ is free as a $\Lambda_0$-module (see \cite[Lemma 2.6]{sz1}). 
We give a proof for the sake of completeness. 

Let $z_j=e_j(u_1^2+u_1^{-2},\ldots, u_m^2+u^{-2}_m)$, ~$x_p=e_p(u_1^2+u_1^{-2},\ldots, u_s^2+u_s^{-2})$, and $ 
y_q=e_q(v_1^2+v_1^{-2},\ldots,v_t^2+v_t^{-2})$.   Then $\mathbb Z[z_1,\ldots, 
z_m]=\mathbb Z[\Lambda_1,\ldots, \Lambda_m]$.  Indeed, 
since $\Lambda_1,\ldots, \Lambda_m$ are expressible as symmetric polynomials in $u_j^2+u_j^{-2}, 1\le j\le m,$ they are expressible as polynomials in 
$z_1,\ldots ,z_m$.  Conversely, since $z_1,\ldots, z_m\in \mathbb  Z[u_1^2, u_1^{-2},\ldots, u_m^2,u_m^{-2}]$ are invariant under  the permutations of the variables $u_1^2,\ldots, u_n^2$ as well as the involutions $u_j^2\mapsto u_j^{-2}$ for every $j$, we see that the $z_j$ belong to the subring of $\mathbb Z[u_1^2, u_1^{-2}, \ldots, u_n^2, u_n^{-2}]$ fixed by the group $\mathbb Z_2^n\rtimes S_n$.   This fixed subring equals $\mathbb Z[\Lambda_1,\ldots, 
\Lambda_m]$; 
see \cite[\S10, Ch. 13]{husemoller}.   So each $z_j$ is expressible as a polynomial in the $\Lambda_i$.

The same argument shows that $\mathbb Z[\lambda_1,\ldots, \lambda_s]=\mathbb Z[x_1,\ldots, x_s]$ and 
$\mathbb Z[\mu_1,\ldots, \mu_t]=\mathbb Z[y_1,\ldots,y_t]$. Consequently, $R=\mathbb Z[\lambda_p,\mu_q;1\le p\le s,1\le q\le t]\subset RH_{n,k}^0.$ 
 
Now using Equation (6) we obtain 
\[ \rho_0(z_j)=\sum_{p+q=j}x_py_q+2\sum_{p+q=j-1}x_py_q, 1\le j\le m-1,\eqno(10)\] 
and $\rho_0(z_m)=2x_s.y_t$ where it is understood that $z_0=x_0=y_0=1$. 
Set $z_1':=z_1-2$, and, inductively, 
$z_r':=z_r-2z_{r-1}', 2\le r<m$, so that $\rho_0(z_r')=\sum_{p+q=r}x_py_q, 1\le r\le m-1$.   
Then $\mathbb Z[z_1',\ldots,z_{m-1}']=
\mathbb Z[z_1,\ldots, z_{m-1}]=\Lambda_0$.  Moreover, we have
\[\rho_0(z_j')=\sum_{p+q=j} x_p.y_q , 1\le j\le m-1.\eqno(11)  \]

The proof that $R$ is a free $\Lambda_0$-module of rank $m-1\choose s$  is now completed using some well-known facts concerning the cohomology of classifying spaces $B\mathrm{U}(s)$ of the unitary group $\mathrm{U}(s)$, as we shall now explain.   We regard 
$R=\mathbb Z[x_1,\ldots, x_s, y_1,\ldots, y_t]$ as a {\it graded} ring where $|x_p|=2p, |y_q|=2q$. Then $\Lambda_0=\mathbb Z[z_1',\ldots, z_{m-1}']$ is a graded subring where $|z_r'|=2r$.   We may identify 
$R$ with $H^*(B(\mathrm{U}(s)\times \mathrm{U}(t));\mathbb Z)$ and $\Lambda_0$ with $H^*(B\mathrm{U}(s+t);\mathbb Z)$ so that 
the inclusion $\Lambda_0\hookrightarrow R$ corresponds to the homomorphism induced by the the projection 
of the fibre bundle $B(\mathrm{U}(s)\times \mathrm{U}(t))\to B\mathrm{U}(s+t)$ with fibre the complex Grassmann manifold $\mathbb CG_{s+t,s}
=U(s+t)/U(s)\times U(t)$. 
The Grassmann manifold bundle is totally non-cohomologous to zero (with $\mathbb Z$-coefficients) and so 
by the Leray-Hirsch theorem $H^*(B(\mathrm{U}(s)\times \mathrm{U}(t));\mathbb Z)$ is a free $H^*(B\mathrm{U}(s+t);\mathbb Z)$-module of rank 
equal to $\textrm{rank}(H^*(\mathbb CG_{s+t,t};\mathbb Z))={s+t\choose s}$.

Since $RH_{n,k}$ is a free $R[\theta]$-module (with basis $\{1,\Delta_{s,t}\}$) by Proposition
 \ref{rhnk},  the last assertion of the lemma follows. 
\end{proof}

\begin{remark}\label{basisB0}
(i) 
We shall denote by $\mathcal B_0$ a basis of $R=\mathbb  Z[\lambda_p,\mu_q;1\le p\le s,1\le q\le t]$ over 
$\Lambda_0$ and assume that $1\in \mathcal B_0$.  Then a $\mathbb Z[\Lambda_1,\ldots,
\Lambda_{m-1}]$-basis for $RH_{n,k}$ is $\mathcal B_0\cup \mathcal B_0\theta\cup \mathcal B_0\Delta_{s,t}
\cup \mathcal B_0\theta\Delta_{s,t}$.  

(ii) The argument in the last paragraph of the above proof is valid irrespective of the parity of $m=s+t+1$.  So 
$R=\mathbb Z[x_1,\ldots,x_s,y_1,\ldots, y_t]$ is a free $\Lambda_0=\mathbb Z[z_1',\ldots, z_{s+t}']$-module for any $s,t\ge 1$.   Moreover,
the quotient ring $R/I$, being isomorphic to $H^*(\mathbb CG_{s+t,s};\mathbb Z)$, is a free abelian group 
of rank ${s+t\choose s}$ where $I$ is the ideal $\langle z_1',\ldots,z_{s+t}'\rangle\subset R$.
\end{remark}

Next we note that irrespective of whether $n\equiv 0$ or $4 \Mod{8},$ we have $\rho((\Delta^+_m)^2-(\Delta^-_m)^2)=0$ and $\rho(\Delta^+_m\Delta^-_m)=\theta\Delta^2_{s,t}=\theta f_{s,t}$.
We have the following consequence of Lemma \ref{freenessoverLambda}.

\begin{lemma} \label{freeness}
The elements $\Lambda_1',\ldots, \Lambda_{m-2}', \rho(\Delta^+_m)\in RH_{n,k}$ are algebraically independent. 
As a module over  $\Lambda:=\mathbb Z[\Lambda_1,\ldots, \Lambda_{m-2},\Delta^+_m]\subset R\Spin(n)$, $RH_{n,k}$ 
is free of rank $2{m-1\choose s}$ with basis $\mathcal B_0\cup \mathcal B_0\theta$.
\end{lemma}
\begin{proof}
Since $\rho(\Delta^+_m)^2=\Delta_{s,t}^2=f_{s,t}$, it suffices to show that $\Lambda_1',\ldots, \Lambda_{m-2}', f_{s,t}$
are algebraically independent in $RH_{n,k}$.   Note that $\Delta^+_m\cdot \Delta^-_m=\Lambda_{m-1}+
\Lambda_{m-3}+\ldots+\Lambda_1$  in $R\Spin(n)$; see \cite[Theorem 10.3, Chapter 14]{husemoller}.  So \[f_{s,t}=\theta\rho(\Delta^+_m\cdot\Delta^-_m)=\Lambda'_{m-1}+\Lambda'_{m-3}+\cdots+1. \]   
Since $\Lambda_1',\ldots , \Lambda_{m-1}'$ are algebraically independent, it follows that $\Lambda_1',
\ldots, \Lambda_{m-2}',f_{s,t}$ are also algebraically independent.  Moreover, we have  $\Lambda'[\rho(\Delta^+_m)]=\rho(\Lambda)\cong \Lambda$.

Let $\mathcal B$ be a basis for $R'[\theta]=R[\theta]$
over $\Lambda'=\mathbb Z[\Lambda_1',\ldots, \Lambda'_{m-1}]$. 
Note that we may take $\mathcal B$ to be $\mathcal B_0\cup \mathcal B_0\theta$ by Remark \ref{basisB0}.
Then 
$\mathcal B$ is a basis for $R[\theta, \rho(\Delta_m^+)]=RH_{n,k}$ 
over $\Lambda'[\rho(\Delta^+_m)]\cong \Lambda$.  In view of Lemma \ref{freenessoverLambda}, 
we conclude that $RH_{n,k}$ is a free module over $\Lambda$ of rank $2{m-1\choose s}$.
\end{proof}

Let $\delta_m=\Delta^+_m-\Delta_m^-$.  Then $R\Spin(n)=\Lambda[\delta_m]$ with $\Lambda$ as in Lemma  \ref{freeness}.  Note that 
$\rho((\Delta^+_m)^2-(\Delta^-_m)^2)=0$ and $\rho(\Delta^+_m\cdot\Delta_m^-)=\theta \Delta_{s,t}^2=\theta f_{s,t}$.  So  the following equations hold in $RH_{n,k}$:
 \[\rho((\Delta^+_m)^2)=\rho(\delta_m^2-2\Delta_m^+\delta_m)=0, \textrm{~and~}\rho(\Delta^+_m)\rho(\delta_m)+(\theta-1)\cdot f_{s,t}=0. \eqno(12)\]

\subsection{Computation of $\mathrm{Tor}^*_{R\Spin(n)}(RH_{n,k},\mathbb Z)$}
We shall apply the change of rings spectral sequence (\S\ref{changeofrings}) to compute 
 $\mathrm{Tor}^*_{R\Spin(n)}(RH_{n,k},\mathbb Z)$.  
In the notation of Theorem \ref{changeofrings},
we let $\Gamma=R\Spin(n),  A=RH_{n,k}, K=C=\mathbb Z$ and $\Lambda
=\mathbb Z[\Lambda_1,\ldots, \Lambda_{m-2},\Delta^+_m]\subset \Gamma=R\Spin(n).$ 
Then $A$ is a free 
$\Lambda$-module via the restriction homomorphism, in view of Lemma \ref{freeness}.    Hence setting 
\[B:=\textrm{Tor}^\Lambda_* (RH_{n,k} ,\mathbb Z),\]
we have, with $\varepsilon\in \{0,1\} $ as in Proposition \ref{restriction}, 
 \[B_q=\textrm{Tor}^\Lambda_q(RH_{n,k},\mathbb Z)=\begin{cases} 
RH_{n,k}/\langle \Lambda_j'-{n\choose j}, 1\le j\le m-2; \theta^\varepsilon 
\Delta_{s,t} -2^{m-1} \rangle, &\textrm{if~}q=0,\\
0,\textrm{~if~} q\ne 0.
\end{cases}\eqno(13)\]  
 Thus 
\[B=B_0=RH_{n,k}/\langle\Lambda_j'-{n\choose j}, 1\le j\le m-2; ~\theta^\varepsilon 
\Delta_{s,t} -2^{m-1} \rangle .\] 

Recall the basis $\mathcal B=\mathcal B_0\cup \mathcal B_0\theta$ of $RH_{n,k}$ over $\Lambda$ given 
in Lemma \ref{freeness}. (See 
Remark \ref{basisB0} for the definition $\mathcal B_0$.) Under the natural projection $\eta: RH_{n,k}\to B$, 
the subring $\rho(\Lambda)$ maps to $\mathbb Z$ and $\mathcal B$ to a $\mathbb Z$-basis $\overline{\mathcal B}=\overline{ \mathcal B}_0\cup \overline{ \mathcal B}_0\theta$ where $\overline{\mathcal B}_0=\eta(\mathcal B_0)$.   It is readily seen that 
$|\overline{\mathcal B}|=|\mathcal B|$. 
We summarise this observation as a lemma.

\begin{lemma}  \label{basisforB} 
The set $\overline{\mathcal B}$ is a $\mathbb Z$-basis for $B$. Thus $B$ is free abelian of rank 
$2{m-1\choose s}$.
\end{lemma}

 By Theorem \ref{changeofrings}, the change of rings spectral sequence collapse 
and we have $\textrm{Tor}^\Gamma_q(A,\mathbb Z)\cong \textrm{Tor}^\Omega_q(B,\mathbb Z)$, where 
$\Omega=R\Spin(n)/\langle \Lambda_j-{n\choose j}, \Delta^+_m-2^{m-1}\rangle =\mathbb Z[\delta_m]$ and 
$\delta_m=\Delta^+_m-\Delta^-_m$.

 Since $\Omega$ is a polynomial ring, 
one can use the Koszul resolution to compute $\textrm {Tor}^\Omega_q(B,\mathbb Z)$.  The $\Omega$-module structure on $B$ is obtained via the algebra homomorphism $\bar\rho:\Omega\to B$ defined by $\rho: R\Spin(n)\to RH_{n,k}$.  In view 
of Proposition \ref{restriction}, we have $\bar\rho(\delta_m)=\epsilon'( \theta-1)\Delta_{s,t}$, where the value of $\epsilon'\in \{1,-1\}$ depends on 
the value of $n$ modulo $8$.   The Koszul resolution of $\mathbb Z$ is 
\[0\to \Omega\cdot\delta\stackrel{d}{\to }\Omega \stackrel {\varepsilon}{\to} \mathbb Z\to 0. \] 
Here $\varepsilon$ is the augmentation defined by $\varepsilon(\delta_m)=0$
and $d(\delta)=\delta_m$.  
Tensoring with the $\Omega$-module $B$ we obtain the following chain complex whose 
homology is $\tor^\Omega_*(B,\mathbb Z)$: 
\[ 0\to B\delta\stackrel{\bar d}{\to} B\to 0\]
where $\bar d(\delta)=\bar d(1\cdot\delta)=\bar \rho(\delta_m)= \epsilon' (\theta-1)\Delta_{s,t}\in B$.
 In particular $\textrm{Tor}^\Omega_q(B,\mathbb Z)=0$ if $q\ge 2$, 
$\textrm{Tor}^\Omega_1(B,\mathbb Z)=\ker(\bar d), \textrm{Tor}^\Omega_0(B,\mathbb Z)=B/\langle (\theta-1)\Delta_{s,t}\rangle$. 

We set 
\[\bar B:=\textrm{Tor}^\Omega_0(B,\mathbb Z)=B/\langle (\theta-1)\Delta_{s,t}\rangle.\eqno(14)\]

Recall from Equation (8) that $\Lambda_j'=\theta^jf_j$ where $f_j=\sum_{0\le p\le j}  \lambda_p\mu_{j-p}\in RH_{n,k}, 1\le j\le m-1, $ and 
\[\lambda_p=\lambda_{k-p}\mathrm{~and~}  \mu_q=\mu_{n-k-q} \mathrm{~when~} p>s, q>t.\]     
Denote by $\eta:RH_{n,k}\to B$ the canonical quotient map and by 
$\bar \eta:RH_{n,k}\to \bar{B}$ the composition 
$RH_{n,k}\stackrel{\eta}{\to }B\to \bar{ B}$ where $B\to \bar{B}$ is the canonical quotient map.  If $x\in RH$, we shall denote $\eta(x)\in B$ by the same symbol  $x$ and we shall denote $\bar\eta(x)\in \bar B$ by $[x]$.

\begin{lemma}  \label{barB} We keep the above notations.   
The following relations hold in $\bar{ B}$: \\
(a) $2^{m-1}([\theta]-1)=0$, $[\Delta_{s,t}]=2^{m-1}$, \\
(b) $\sum_{0\le p\le j }[\lambda_{p}][\mu_{j-p}]=[f_j]= {n\choose j}[\theta^j], 1\le k\le m-1,$  (where $[\lambda_p]=[\lambda_{k-p}], [\mu_q]=[\mu_{n-k-q}]$),\\
(c) $[\Delta_{s,t}^2]=(\sum_{0\le p\le s}[\lambda_p])(\sum_{0\le q\le t}[\mu_q])=[f_{s,t}]=2^{2m-2}$.\\
\end{lemma}
\begin{proof}   (a).   We have, by Proposition \ref{restriction}, $\rho(\Delta^+_m)=\theta^\varepsilon \Delta_{s,t}, $ in $RH_{n,k}$ where $\varepsilon\in \{0,1\}$ depending on the value of $n$ modulo $8$. 
Since $([\theta]-1)[\Delta_{s,t}]=0$ in $\bar B$, irrespective of the value of $\varepsilon$ we have $\bar\eta \circ \rho(\Delta_m^+)=[\Delta_{s,t}]$ in $\bar B$.  On the other hand, 
since $\Delta^+_m=2^{m-1}$ in $\Omega$, we obtain that $2^{m-1}=\eta\rho(\Delta^+_m)=\theta^\varepsilon \Delta_{s,t}$ in $B$.    It follows that $[\Delta_{s,t}]=2^{m-1}$ and so $2^{m-1}([\theta]-1)=0$.

(b).   It is clear that, when $1\le j\le m-2$, the 
relation $f_j=\bar\rho(\Lambda_j)\theta^j={n\choose j}\theta^j$ holds in $B$ and hence in $\bar B$ using $\theta^2=1$. 
Since $\Delta_m^+\Delta_m^-=\sum_{1\le j\le m}\Lambda_{2j-1}$ in $R\Spin(n)$, and since 
$\bar\eta\circ \rho(\Delta^\pm_m)=[\Delta_{s,t}]=[\theta][\Delta_{s,t}]=2^{m-1} $ in $\bar B$,  applying $\bar\eta \circ \rho$ we obtain the following equations in $\bar B$:
\[
\begin{array}{rcl} 2^{2(m-1)} &=& \bar\eta \circ \rho(\Delta^+_m\Delta^-_m)\\
&=&\bar\eta\circ  \rho(\sum_{1\le j\le m} \Lambda_{2j-1})\\
&=&[f_{m-1}]-{2m\choose m-1}+\sum_{1\le j<m/2} {2m\choose 2j-1}\\
&=& [f_{m-1}]-{2m\choose 2m-1} +2^{2m}/4\\
\end{array}
\]
since 
$\sum_{1\le j<m/2}{2m\choose 2j-1}=(1/2)\sum_{1\le j\le m}{2m\choose 2j-1}=2^{2m}/4$.  Hence 
$[f_{m-1}]={2m\choose m-1}$.

(c).   
Since $\Delta_{s,t}^2=f_{s,t}$ holds in $B$, and since $[\Delta_{s,t}]=2^{m-1}$ holds in $\bar B$, we see 
that $[f_{s,t}]=2^{2m-2}$ in $\bar B$.
\end{proof}

\begin{remark}
It turns out that the relation (c) is a consequence of relations (a), (b).  Indeed, 
recalling that $[\lambda_p]=[\lambda_{ k-p}], [\mu_q]=[\mu_{n-k-q}]$ in $\bar B$ and also that $k=2s+1, n-k=2t+1$, we have 
\[ \begin{array}{rcl}
f_{s,t}=[\Delta_{s,t}^2]&=&(\sum_{0\le p\le s}[\lambda_p])(\sum_{0\le q\le t}[\mu_q] )\\
&=&(1/4)(\sum_{0\le p\le k}[\lambda_p])(\sum_{0\le q\le n-k}[\mu_q])\\
&=&(1/4) \sum_{0\le r\le n}(\sum_{0\le j\le r}[ \lambda_j][\mu_{r-j}])\\
&=&(1/4)\sum_{0\le r\le n}{n\choose r}[\theta]^r, ~\mathrm{~using (b)},\\
&=&(1/4)(1+[\theta])^n.\\
\end{array}
\]
Since $[\theta]^2=1$, we have $(1+[\theta])^2=2(1+[\theta])$.  So $(1+[\theta])^{r}=2^{r-1}(1+[\theta])$ for all $r\ge 1$.   Therefore, since $n=2m\ge 4 $, 
we have 
\[\begin{array}{rcl}
(1/4)(1+[\theta])^n&=&(1/4)(1+[\theta])^3\cdot(1+[\theta])^{n-3}\\
&=&(1+[\theta])\cdot (1+[\theta])^{2m-3}\\
&=&(1+[\theta])^{2m-2}\\
&=&2^{2m-3}(1+[\theta])\\
&=& 2^{2m-2},
\end{array}
\] using $2^{2m-3}[\theta]
=2^{2m-3}$.     Therefore $f_{s,t}=2^{2m-2}$.  
\end{remark}

\begin{lemma}\label{rankandtorisionofbarB}
With the above notations, the rank of the abelian group $\bar B$ equals ${m-1\choose s}$.  Moreover 
the torsion subgroup of $\bar B$ is generated as a $B$-module by $(\theta-1)$.  In particular, any 
torsion element has order $2^r$ for some $r\le m-1$.  
\end{lemma}
\begin{proof} 
In view of Lemma \ref{basisforB}, the set $\overline{\mathcal B}_0\cup \overline{\mathcal B}_0(\theta-1) $ is a basis 
for $B$. Under the quotient map $B\to \bar B$, the abelian group $\bar B_0$ generated by $\overline{\mathcal B}_0$ projects 
isomorphically onto a summand of $\bar B_0$.  Since $2^{m-1}([\theta]-1)=0,$ 
the subgroup $C$ of $\bar B$ is generated by 
$([\theta]-1)\overline{\mathcal B}_0$ 
consists only of elements whose (additive) order divides $2^{m-1}$.  This completes the proof. 
\end{proof}

We now turn to $\textrm{Tor}^\Omega_1(B,\mathbb Z)=\ker (\bar d:B\delta\to B)$.  Since $\bar d(\delta)=\pm(\theta-1)
 \Delta_{s,t}$, $\ker(\bar d)$ is the $B$-submodule $J\cdot\delta$ where $J\subset B$ is the annihilator ideal of 
$(\theta-1)\Delta_{s,t}\in B.$  
 It is clear that $(\theta+1)\in J$ since $\theta^2-1=0$.   We claim that  $J$ {\it equals} the ideal generated by 
$\theta+1$.
In order to see this, let $x\in J$ and let 
$\bar {\mathcal B}_0=\{b_j\}$.   Write 
$x=\sum y_jb_j+\theta \sum  z_jb_j$ where $y_j,z_j\in \mathbb Z$.  
Since $x\in J$, multiplying by $(\theta-1) \Delta_{s,t}$,  and using the relations $\Delta_{s,t}=2^{m-1}\theta^\varepsilon$ 
(where the value of $\varepsilon\in \{0,1\}$ depends on the parity of $m$) 
and $\theta(\theta-1)=1-\theta$ in $B$, we obtain that 
\[ 2^{m-1} (\theta-1)\theta^\varepsilon \sum y_j b_j  -2^{m-1}\theta^{\varepsilon} (\theta-1)\sum z_jb_j=0.\]   
Since $B$ is a free abelian group, and since $\theta^\varepsilon$ is invertible in $B$, 
the above equation can be rewritten as 
$-(\sum( y_j-z_j)b_j)+ \theta \sum( y_j-z_j)b_j=0$. This implies that $y_j=z_j$ for all $j$ . Therefore $x=(\theta+1)
(\sum y_jb_j)\in J.$

We are now ready to prove Theorem \ref{main}.

\noindent
{\em Proof of Theorem \ref{main}}:
The Hodgkin spectral sequence $\mathrm{Tor}^*_{R\Spin(n) }(RH_{n,k},\mathbb Z)$ converges to $K^*(G_{n,k})$.
Since $\mathrm{Tor}_{R\Spin(n)}^*(RH_{n,k},\mathbb Z)\cong \mathrm{Tor}^* _\Omega(B,\mathbb Z)$, and since 
$\tor^*_\Omega(B,\mathbb Z)$ is generated by degree $-1$ elements, by the discussion in \S2 we obtain 
that $K^0(G_{n,k}) =\tor^0_\Omega(B,\mathbb Z)=\bar B$ and $K^{-1}(G_{n,k})=\tor_1^\Omega(B,\mathbb Z)=
\mathrm{Ann}(\theta-1)\subset B$.  The theorem now follows from Equation (14),  Lemma \ref{barB}, and the above discussion 
that describes $\mathrm{Ann}((\theta-1)\Delta_{s,t})$.  \hfill $\Box$

Let $\xi=\xi_{n,k}$ be the Hopf line bundle over $G_{n,k}$.  It is associated to the double cover $\widetilde{G}_{n,k}\to G_{n,k}$.   If $\eta$ is a real vector bundle, we denote by $\eta^\mathbb C$ the complexification of $\eta$.   Note that 
$\eta^\mathbb C$, regarded as a real vector bundle via restriction of 
scalars, is isomorphic to $\eta\oplus \eta$. See \cite[p. 176]{milnor-stasheff}. 
\begin{proposition} 
 Let $n=2m, k=2s+1$. If 
$n\equiv 0\pmod 4, k\equiv 1\Mod{ 2}$ and $k(n-k)<2^m$,  then $2^m\xi\cong 2^m\epsilon_\mathbb R$ where $n=2m$ but $[2^{m-2}\xi]\ne 2^{m-2}$ in $KO(G_{n,k}).$  
If $n\equiv 0\Mod{ 8}$ and $k(n-k)<2^{m-1}$, then $2^{m-1}\xi\cong 2^{m-1}\epsilon_\mathbb R$ provided 
$k(n-k)<2^{m-1}$. 
\end{proposition}
\begin{proof} 
Since $2^{m-1}[\xi^\mathbb C]=2^{m-1}\theta=2^{m-1}\in K(G_{n,k})$, it follows that $2^{m}[\xi]= 2^{m}\in KO(G_{n,k})$.  
If $\dim G_{n,k}=k(n-k)<2^m=\textrm{rank}(2^m\xi)$, then equality of the classes of the vector 
bundles $[2^m\xi]$ and $[2^m\epsilon_\mathbb R]=2^m$ in $KO(G_{n,k})$ implies  the {\it isomorphism} of the 
vector bundles: $2^m\xi\cong 2^m\epsilon_\mathbb R$.   See \cite[Theorem 1.5, Chapter 8]{husemoller}. 

When $n\equiv 0\Mod{8}$, the representations $\Delta^+_m,\Delta_m^-\in R\Spin(n)$ are real, that is, they arise 
as complexification of {\it real }representations $\Delta^+_{m,\mathbb R}, \Delta^-_{m,\mathbb R}$ of $\Spin(n)$. 
See \cite[\S12, Chapter 13]{husemoller}.  Evidently $\theta$ is real.  Indeed 
$\theta=\chi \otimes _\mathbb R\mathbb C$ of $H_{n,k}$ where $\chi:H_{n,k}\to \O(1)$ is defined by 
the projection $H_{n,k}\to H_{n,k}/H^0_{n,k}\cong \O(1)$.   

The line bundle associated to $\chi$ is isomorphic to $\xi$ whereas the bundle associated 
to $\Delta_{m,\mathbb R}^-$ equals the trivial real vector bundle of rank $2^{m-1}$. 
This can be shown to imply that $2^{m-1}[\xi]=2^{m-1}\in KO(G_{n,k})$.  As before, this leads 
to the isomorphism $2^{m-1}\xi\cong 2^{m-1}\epsilon_\mathbb R$ when $k(n-k)<2^{m-1}$.   
\end{proof}

As for the torsion part of $K^0(G_{n,k})$, it has no $p$-torsion for any odd prime $p$.   
For any $n,k$, the   element 
$[\Lambda^k(\gamma_{n,k}^\mathbb C)]-1=[\xi^\mathbb C]-1\in K(G_{n,k})$ generates a finite cyclic subgroup of order $2^r$ for some $r$.  There are the obvious inclusions $i:G_{n,k}\hookrightarrow G_{n+1,k+1}, j: G_{n,k}\hookrightarrow G_{n+1,k}$ which have the property that 
$i^*(\gamma_{n+1,k+1})\cong \gamma_{n,k}\oplus \epsilon_\mathbb R$ and  $j^*(\gamma_{n+1,k})=\gamma_{n,k}. $

\begin{theorem} \label{hopforder} Suppose that $n=4l+j, k=2s+\varepsilon, 1\le j\le 3, \varepsilon\in\{0,1\}$.   Let $2^r$ be the order of $[\xi^ \mathbb C]\in K(G_{n,k})$.   Then $2l-1\le r\le 2l+1$. 
\end{theorem}
\begin{proof}
Suppose $\varepsilon=1$.  Then we have inclusions $G_{4l,k}\stackrel{j_0}{\hookrightarrow} G_{4l+j,k}\stackrel{j_1}{\hookrightarrow }G_{4l+4,k}$ where $j_1^*(\xi_{4l+4,k})=\xi_{n,k},~ j_0^*(\xi_{n,k})=\xi_{4l,k}$.
The bounds for $r$ now follow from Theorem \ref{main}. 

When $\varepsilon=0$, we use the inclusions $G_{4l,2s-1}\stackrel{i_0}{\hookrightarrow}  G_{n,k}\stackrel{i_1}{\hookrightarrow}G_{4l+4,2s+1}$.  
When $s=1$, $G_{4l,2s-1}=\mathbb RP^{4l-1}$ and the order of the bundle $[\xi^ \mathbb C]-1$ 
is known to be $2^{2l-1}$ from the work of Adams \cite[Theorem 7.3]{adams}.   Now we proceed 
exactly as in the case $\varepsilon=1$. 
\end{proof}

\section{$K$-theory of $G_{n,k}$ for arbitrary values of $n,k$}\label{nkarbitrary}
In this section we shall prove Theorem \ref{main2}. The proof will make use of the Chern character $\textrm{ch}:K^*(G_{n,k})\otimes \mathbb Q\to 
H^*(G_{n,k};\mathbb Q)$.  
We begin by recalling, in Theorem \ref{cohomologyalgebra} and the following paragraph, the rational cohomology algebra of 
 the Grassmann manifolds.   We refer the reader to \cite[\S 15]{milnor-stasheff} for the definition and 
properties of Pontrjagin classes.   We shall write $k=2s+\varepsilon, n-k=2t+\eta$ where $\varepsilon,\eta\in \{0,1\}$ so that $n=2s+2t+\varepsilon+\eta$.

We denote by $\beta_{n,k}$ the canonical $(n-k)$-plane bundle over $G_{n,k}$ whose fibre over $L\in G_{n,k}$ is the vector space $L^\perp\subset \mathbb R^n$. 
We have 
$ \gamma_{n,k}\oplus \beta_{n,k}\cong n\epsilon_\mathbb R,$ and, (denoting the 
complexification $\gamma_{n,k}\otimes \mathbb C$ by $\gamma_{n,k}^\mathbb C$ etc.,) we obtain  
\[\gamma_{n,k}^\mathbb C\oplus \beta_{n,k}^\mathbb C=n\epsilon_\mathbb C. \eqno(15)\]

Let $p_j=p_j(\gamma_{n,k})\in H^{4j}(G_{n,k};\mathbb Z[1/2])$, $1\le j\le s$, be the $j$th (rational) Pontrjagin class of $\gamma_{n,k}$, and let $q_j=p_j(\beta_{n,k}), 1\le j\le t$.  Since $\gamma_{n,k}\oplus \beta_{n,k}\cong 
n\epsilon_\mathbb R$, we have, for $1\le r\le s+t,$
\[ \sum_{0\le j\le s}  p_jq_{r-j}=0,\eqno(16)\]
where it is understood that $p_0=q_0=1, p_i=0, q_j=0$ if $i>s,j>t$. 
In fact, the cohomology algebra $H^*(G_{n,k};\mathbb Z[1/2])$ has the following description.   It can be derived from the known description of $H^*(\widetilde{G}_{n,k};\mathbb Z[1/2])$ as the fixed subring under the action of the deck transformation group of the double covering $\widetilde{G}_{n,k}\to G_{n,k}$.   We refer the reader to \cite[Theorem 15.9]{milnor-stasheff} for the description of $H^*(\widetilde{G}_{n,k};\mathbb Z[1/2])$.

\begin{theorem}  \label{cohomologyalgebra} With the above notations, we have
\[H^*(G_{n,k};\mathbb Z[1/2])=\mathbb Z[1/2][p_1,\ldots, p_s;q_1,\ldots, q_t, v_{n-1}]/J \eqno(17) \]
where degree of $v_{n-1}=n-1$, and the ideal $J$ is generated by the following elements:\\
(i) $\sum_{0\le j\le r} p_jq_{r-j}, 1\le r\le s+t$,\\
(ii) $v_{n-1}$ if $n$ is odd or $k$ is even; $v_{n-1}^2$ if $n$ is even and $k$ odd.
\end{theorem}

As a consequence we note that $H^*(G_{n,k};\mathbb Z)$ has no $p$-torsion except when $p=2$.  

Denote by $P_{n,k}\subset H^*(G_{n,k};\mathbb Q)$ the even-graded subalgebra, namely, 
$H^{\textrm{ev}}(G_{n,k};\mathbb Q)=\oplus_{r\ge 0}H^{2r}(G_{n,k};\mathbb Q)=
\mathbb Q[p_1,\ldots,p_s;q_1,\ldots, q_t]/\!\!\sim$.  Then  $P_{n,k}$ depends only on $s,t$ and not on the values of $\varepsilon,\eta\in \{0,1\}$, and,  $\dim_\mathbb Q
P_{n,k}={s+t\choose s}$. 
Moreover, $P_{n,k}=H^*(G_{n,k};\mathbb Q)$, except when $n=2s+2t+2$ is even and $k=2s+1$ is odd.   When $n=2s+2t+2, k=2s+1,$ we have $H^{\textrm{odd}}(G_{n,k};\mathbb Q)=v_{n-1}P_{n,k}
\cong P_{n,k}$ as a $P_{n,k}$-module.   We have a natural $\mathbb Z_2$-gradation on $H^*(G_{n,k};\mathbb Q)$ 
defined by the parity of the degree.  

Recall the Chern character map $\textrm{ch}:K^*(G_{n,k})\otimes \mathbb Q\to H^*(G_{n,k};\mathbb Q)$, which is 
an isomorphism of $\mathbb Z_2$-graded rings.   
So $K^0(G_{n,k})$ has rank equal to $\dim_\mathbb Q P_{n,k}={s+t\choose s}$.   
In case $n$ is odd or $k$ is even, we have $H^{\textrm{odd}}(G_{n,k};\mathbb Q)=0$ and so 
$K^1(G_{n,k})$ 
is a finite abelian group.   When $n$ is even and $k$ is odd, $K^1(G_{n,k})$ has rank equal to that of $K^0(G_{n,k})$.  

We now turn to the proof of Theorem \ref{main2}.   We shall denote by $\phi$ the inclusion map 
$\mathcal K_{n,k}\hookrightarrow K(G_{n,k})$.   

\noindent
{\it Proof of Theorem \ref{main2}}   The inclusion $\phi:\mathcal K_{n,k}\hookrightarrow K(G_{n,k})$ induces an  inclusion $\phi\otimes 1:\mathcal K_{n,k}\otimes \mathbb Q\to K(G_{n,k})\otimes \mathbb Q$.

 We need to show that 
the composition $\textrm{ch}\circ (\phi\otimes 1): \mathcal K_{n,k}\otimes \mathbb Q\to P_{n,k}$  is surjective. 
Note that, in view of Equation (16), $P_{n,k}$ is generated by $p_j,1\le j\le s$.
So we need only show that the $p_j\in P_{n,k}$ are in the image of $\textrm{ch}\circ (\phi\otimes 1)$. 

We have a formal expression of $p_j=p_j(\gamma_{n,k})$ in terms of 
the Chern `roots'  $x_j, -x_j, 1\le j\le s,$ of $\gamma_{n,k}^\mathbb C$ given as 
$p_j=(-1)^je_{j}(x_1^2,\ldots, x_k^2), 1\le j\le s,$ where $e_j$ denotes the $j$th elementary symmetric polynomial in the indicated arguments.  (See \cite[\S15]{milnor-stasheff}.)   
From the definition of Chern character we have 
$\textrm{ch}(\gamma_{n,k}^\mathbb C)=k+2\sum_{m\ge 1} \sum_{1\le j\le s} x_j^{2m}/(2m)!
=k+2\sum_{m\ge 1} u_m/(2m)!$ where $u_m:=\sum _{1\le m\le s} x_j^{2m}$ for $m\ge 1$.    
The 
symmetric polynomials can be expressed as polynomials in the power sums over $\mathbb Q$ and so we have 
\[(-1)^jp_j=u_j/j+F_j(u_1, \ldots, u_{j-1}), ~~1\le j\le s,\eqno(18)\] where $u_0=k$
and $F_j(u_1,\ldots, u_{j-1})\in H^{4j}(G_{n,k};\mathbb Q)$ is a suitable polynomial in $u_1,\ldots, u_{j-1}$.  
So it suffices to show that the $u_j$ are in the image of $\textrm{ch}\circ (\phi\otimes 1).$   
To see this, it is convenient to use the Adams operations $\psi^r$.  Note that $\mathcal K_{n,k}$ contains 
$\Lambda_j(\gamma_{n,k}^\mathbb C)$ 
and so it also contains $\psi^r(\gamma_{n,k}^\mathbb C)$ for all integers $r\ge 1$ since the $\psi^r$ can be expressed 
(with $\mathbb Z$-coefficients) in terms of the exterior power operations.   Although 
$\psi^r(\gamma^\mathbb C_{n,k})$ is only a virtual bundle, its Chern characters are easy to compute   
since $rx_j, -rx_j$ are its Chern roots. Thus, writing $d=\lfloor k(n-k)/2\rfloor $, 
we have, for $r\in \mathbb Z$, 
 \[\begin{array}{rcl}
 v_r&:=&\textrm{ch}([\psi^r(\gamma_{n,k}^\mathbb C)] -k)\\
 &=&2\sum_{m\ge 1}(\sum_{1\le j\le s} r^{2m} x_j^{2m}/(2m)!)\\
&=&2\sum_{1\le m\le d} r^{2m}u_m /(2m)!.\\
\end{array}
\eqno(19)\]

We obtain the equation 
$2uM=v$ where $M=(m_{ij})$ is the $d\times d$ 
 matrix defined as $m_{ij}=j^{2i}$, and $u=(u_1/2!,u_2/4!,\ldots, u_d/(2d)!), v=(v_1,\ldots, v_d)$ are regarded as (row) vectors in the $d$-fold direct sum $(H^{\textrm{ev}}(G_{n,k};\mathbb Q))^d$.  Since $M$ is invertible and since the 
$v_j$ are in the image of $\textrm{ch}\circ( \phi\otimes 1)$, it follows that the $u_j/(2j)!$ are also in the image of 
$\textrm{ch}
\circ (\phi\otimes 1)$ for $1\le j\le d$.   So $u_1,\ldots, u_s$ are in the image of $\textrm{ch}\circ( \phi\otimes1)$.  This completes the proof.   \hfill $\Box$

We conclude by giving, in Proposition \ref{almostK}, a description of $\mathcal K_{n,k}$ as a quotient of a ring $K_{n,k}$, explicitly described 
in terms of generators and relations, with finite kernel.  It seems plausible that $K_{n,k}$ is isomorphic to $\mathcal K_{n,k}$ but we have not been able to prove this.


The operator $\Lambda_t=\sum_{r\ge 0} \Lambda^r t^r$, which is a formal power series in the indeterminate $t$ whose coefficients are exterior power operators, has the property $\Lambda_t(\omega_0\oplus\omega_1)=\Lambda_t(\omega_0).\Lambda_t(\omega_1)$ for any two complex vector bundles $\omega_0,\omega_1$.  So we have $\Lambda_t(\gamma^\mathbb C_{n,k})\cdot \Lambda_t(\beta_{n,k}^\mathbb C)=(1+t)^n$ since $\Lambda_t(\epsilon_\mathbb C)=(1+t)$.  Equivalently, for any $r\ge 1$, we have 
\[
\sum_{p+q=r}\Lambda^p(\gamma_{n,k}^\mathbb C)\otimes \Lambda^{q}( \beta_{n,k}^\mathbb C))=
{n\choose r}. \eqno(20)\]

We know that $2^r\xi^\mathbb C$ is stably trivial for some $r$ where $\xi=\xi_{n,k}$ denotes the Hopf line bundle over $G_{n,k}=\SO(n)/S(\O(k)\times \O(n-k))$.   By Theorem \ref{hopforder}, one may take $r=m+1$.  
We let $\nu $ be the least positive integer for which this happens.  
Then $(1-[\xi^\mathbb C])^{\nu+1}=2^\nu(1-[\xi^\mathbb C])=0$ in $K(G_{n,k})$.  Note that 
$\xi=\Lambda^k(\gamma_{n,k})=\Lambda^{n-k}(\beta_{n,k})$ is associated to the
character $\chi:  S(\O(k)\times \O(n-k))\to \O(1)$  defined as $\bigl(\begin{smallmatrix} A&0
\\0&B\end{smallmatrix}\bigr)\mapsto \det(A)$.  We let $\theta$ be the complexification of $\chi$ so that $\xi^\mathbb C$ is associated to $\theta$.  
We shall denote $[\xi^\mathbb C]\in K(G_{n,k})$ by $[\theta]$.

For any real vector space $V$ of dimension $k$, one has a functorial non-degenerate bilinear pairing 
$\Lambda^p(V)\times \Lambda^{k-p}(V)\to \Lambda^k(V)$ defined as $(u,v)\mapsto u\wedge v$.  If $V$ is an inner product space, then we have the induced inner product $\Lambda^q(V)\times \Lambda^q(V)\to \mathbb R$ defined as $(u_1\wedge \cdots \wedge u_q, v_1\wedge 
\cdots \wedge v_q)\mapsto \det ((u_i,v_j))$.   Thus, we obtain 
a natural isomorphism $\Lambda^p(V)\cong \Lambda^{k-p}(V)\otimes \Lambda^k(V)$.   This yields 
an isomorphism $\Lambda^p(\gamma_{n,k})\cong \Lambda^{k-p}(\gamma_{n,k})\otimes \xi$ of real vector bundles. 
See \cite[\S 2]{milnor-stasheff}.   A similar isomorphism holds for $\beta_{n,k} $ as well.  
Complexifying we obtain the following isomorphisms for $1\le p\le k, 1\le q\le n-k$:
 \[
\Lambda^p(\gamma^\mathbb C_{n,k})\cong \xi^\mathbb C\otimes\Lambda^{k-p}(\gamma^\mathbb C_{n,k}), ~\Lambda^q(\beta^\mathbb C_{n,k})\cong 
\xi^\mathbb C\otimes \Lambda^{n-k-q}(\beta_{n,k}^\mathbb C).\eqno(21)\]

We are now ready to define the ring $K_{n,k}$.

\begin{definition} \label{defknk}
Let $A=\mathbb Z[\theta]/\langle \theta^2-1, 2^\nu(1-\theta)\rangle$.   Then $A\cong \mathbb Z\oplus \mathbb Z_{2^r}(1-\theta)$.  
Write $k=2s+\varepsilon, n-k=2t+\eta$ where $\varepsilon,\eta\in \{0,1\}$ so that $n=2s+2t+\varepsilon+\eta$.  
We define $K_{n,k}:=A[\lambda_1,\ldots, \lambda_k, \mu_1,\ldots, \mu_{n-k}]/I$,  the quotient of the polynomial algebra over $A$ where the ideal $I$ is generated by the following elements:\\
(i) $\lambda_{k-p}-\theta\lambda_p, ~\mu_{k-q}-\theta\mu_q$ for $1\le p\le k, ~1\le q \le n-k$,  \\
(ii)   $Q_r(\lambda,\mu)-{n\choose r}$ for $1\le r\le n$ where $Q_r(\lambda, \mu):=\sum_{p+q=r, 0\le p \le k,0\le q\le n-k} \lambda_p \mu_{q}$, for 
$1\le r\le n$.
\end{definition} 

\noindent
\begin{remark} \label{observations}
(a) The $A$-algebra $K_{n,k}$ is generated by  $\lambda_1,\ldots, \lambda_s, \mu_1,\ldots, \mu_t$.  This is 
immediate from the relations \ref{defknk}(i).  \\
(b) In fact, using the relations \ref{defknk} (ii), (and (a)), we see that $\mu_1=n-\lambda_1$, and, if $2\le r\le t$, then 
$\mu_r$ can be expressed in terms of the $\lambda_1,\ldots, \lambda_s, \mu_1,\ldots,\mu_{r-1}$ (with coefficients in $A$).   So, by induction, the 
$\mu_r$ can be expressed in terms of $\lambda_1,\ldots, \lambda_s$.  Hence $K_{n,k}$ is generated by $\lambda_p, 
1\le p\le s$.\\
(c) One has a ring homomorphism $A\to \mathbb Z$ which maps $\theta$ to $1$ with kernel the ideal $A(1-\theta).$
\end{remark} 
Set $\bar K_{n,k}:=K_{n,k}\otimes_A\mathbb Z=K_{n,k}/(1-\theta)K_{n,k}=\mathbb Z[\lambda_1,\ldots, \lambda_s,\mu_1,\ldots, \mu_{t}]/I_0$ where $I_0$ is the ideal generated by the elements listed in Definition \ref{defknk} (ii), and where 
$\lambda_p=\lambda_{k-p}, \mu_q=\mu_{n-k-q}$ for $p>s, q>t$ .

\begin{lemma}
One has the following isomorphisms  of rings:
\[\bar{K}_{2s+2t+2, 2s+1}\stackrel{\alpha_0}{\to} \bar{K}_{2s+2t+1, 2s+1}\stackrel{\alpha_1}{\to}\bar K_{2s+2t,2s} ,  \eqno(22)\]
where,
$\alpha_0(\lambda_p)=\lambda_p, ~\alpha_0(\mu_q)=\mu_q+\mu_{q-1},$ and, 
$\alpha_1(\lambda_p)=\lambda_p+\lambda_{p-1},~ \alpha_1(\mu_q)=\mu_q$, for all  $ p\le k, q\le n-k$.  (It is understood that $\lambda_0=1=\mu_0$.)   As an abelian group $\bar K_{n,k}$ is free of rank ${s+t\choose s}$ 
where $(n,k)=(2s+2t+2,2s+1), (2s+2t+1,2s+1),(2s+2t, 2s)$.    

\end{lemma}
\begin{proof}
It is readily verified that 
$\alpha_0,\alpha_1$ are surjective homomorphisms.    We need to show that they are injective as well.  

Consider 
 $\beta_0:\bar K_{2s+2t+1, 2s+1}\to \bar K_{2s+2t+2,2s+1}$,  and,  $\beta_1: \bar K_{2s+2t, 2s}\to \bar K_{2s+2t+1, 2s+1}$ defined as follows:  for $p\le s, q\le t$, \\
$\beta_0(\lambda_p)=\lambda_p, ~
~\beta_0(\mu_q)= \sum_{0\le j\le q}(-1)^{q-j} \mu_j$, and \\
$\beta_1(\lambda_p)= \sum_{0\le j\le p}(-1)^{p-j}\lambda_j,~\beta_1(\mu_q)=\mu_q$. \\ 
Straightforward verification, using the identity $\sum_{0\le j\le r}(-1)^j{n\choose r-j}={n-1\choose r}$, 
shows that $\beta_0$ and $\beta_1$ are well-defined homomorphisms of rings.  Again, 
these are surjective, since the generators $\lambda_p$ (resp. $\mu_q$) are in the image of $\beta_0$ 
(resp. $\beta_1$). 

We claim that $\alpha_0,\beta_0$ (resp. $\alpha_1, \beta_1$) are inverses of each other.  
Indeed, $\beta_0\circ \alpha_0(\lambda_p)=\lambda_p$ for all $p\le s$ and $\alpha_0\circ \beta_0(\lambda_p)=
\lambda_p$ for all $p$. By Remark \ref{observations}(b) above, our claim follows.
Similarly $\alpha_1,\beta_1$ are inverses of each other.

For the last assertion, we need only consider the case $(n,k)=(2s+2t+2,2s+1)$.  
The ring  $\bar K_{2s+2t+2,2s+1}$ is isomorphic to the quotient ring $R/I\cong H^*(\mathbb C G_{s+t,s};\mathbb Z)$ 
considered in Remark \ref{basisB0}(ii).  Hence $\bar K_{2s+2t+2,2s+1}$ is a free abelian group 
of rank ${s+t\choose s}$.
\end{proof}

\begin{proposition}\label{almostK}
One has a surjective homomorphism of rings $\kappa:K_{n,k}\to \mathcal K_{n,k}$ with finite kernel, 
defined as 
$\kappa(\lambda_j)=[\Lambda^j(\gamma_{n,k}^\mathbb C)], 1\le j\le k$.
\end{proposition}
\begin{proof}
In view of Equations (20) and (21), $\kappa$ is a well-defined ring homomorphism.  Clearly $\kappa(\lambda_j)
=[\Lambda^j(\gamma_{n,k}^\mathbb C)]$ for all $j$ and so, by the definition of $\mathcal K_{n,k}$, $\kappa$ 
is surjective.  Since both $K_{n,k}, \mathcal K_{n,k}$ have the same (finite) rank, it follows that $\ker(\kappa)$ is 
finite. 
\end{proof}

\end{document}